\documentclass{article}
\usepackage{amsthm,color}
\usepackage{epsfig}
\usepackage{psfrag,wrapfig}

\def\silent#1\par{}

\silent

\psfrag{y0}{$y_0$}
\psfrag{y1}{$y_1$}
\psfrag{y2}{$y_2$}
\psfrag{yi}{$y_i$}
\psfrag{yj}{$y_j$}
\psfrag{ym}{$y_m$}
\psfrag{x1}{$x_1$}
\psfrag{x2}{$x_2$}
\psfrag{xj}{$x_j$}
\psfrag{xm}{$x_m$}
\psfrag{xi+1}{$x_{i+1}$}
\psfrag{xj+1}{$x_{j+1}$}

\gdef\BalAbra#1#2#3#4{\begin{wrapfigure}{L}{#2mm}\epsfig{file=#1.eps,width=#2mm}
\caption{{\small#3}}\label{#4}
\end{wrapfigure}}

\makeatletter
\renewcommand{\@seccntformat}[1]{\@nameuse{the#1}.\quad}
\makeatother

\def\boxit#1{\medskip\vbox{\hrule \hbox{\vrule\kern15pt\vbox{\kern5pt
\vbox{\advance\hsize -30pt #1\par} \kern5pt}\kern15pt\vrule}\hrule} \medskip}

\def\rema#1//{\medskip {\bf #1}~}

\def\eqref#1{(\ref{#1})}

\def\Case#1){\medskip\noindent{\bf Case #1)}~}

\newtheorem{theorem}{Theorem}[section]
\newtheorem{conjecture}[theorem]{Conjecture}
\newtheorem{corollary}[theorem]{Corollary}
\newtheorem{claim}[theorem]{Claim}

\newtheorem{definition}[theorem]{Definition}
\newtheorem{lemma}[theorem]{Lemma}
\newtheorem{problem}[theorem]{Problem}

\def\Proof#1. {\medskip\noindent{\bf Proof #1.}}

\def\cH{{\mathcal H}}

\begin{document}

\hsize=370pt

\title{\bf On  $3$-uniform hypergraphs without linear cycles\thanks{The authors are grateful for the hospitality of the  Mittag-Leffler Institute program Graphs, Hypergraphs and Computing, during which this research was conducted.}}
\author{
A. Gy\'arf\'as, \thanks{Research supported in part by
the OTKA Grant No. K104343.}\\
\small Alfr\'ed R\'enyi Institute of Mathematics,
 Hungarian Academy of Sciences\\ [-0.8ex]
\small 1053 Re\'altanoda u 13-15, Budapest, Hungary \\
\small\texttt{gyarfas.andras@renyi.mta.hu}
\and E. Gy\H ori, \thanks{Research supported in part by
the OTKA Grant No. K101536.}\\
\small Alfr\'ed R\'enyi Institute of Mathematics, Hungarian Academy of Sciences\\[-0.8ex]
\small 1053 Re\'altanoda u 13-15, Budapest, Hungary\\\small
\texttt{gyori.ervin@renyi.mta.hu}
\and M. Simonovits,\thanks{Research supported in part by
the OTKA Grant No. K101536 and ERC-AdG. 321104.}\\\\[-0.8ex]
\small Alfr\'ed R\'enyi Institute of Mathematics, Hungarian Academy of Sciences\\[-0.8ex]
\small 1053 Re\'altanoda u 13-15, Budapest, Hungary \\\small
\texttt{simonovits.miklos@renyi.mta.hu}}
\maketitle
\begin{abstract}
We explore properties of $3$-uniform hypergraphs $H$ without linear cycles. Our main results are that these hypergraphs must contain a vertex of strong degree at most two and must have independent sets of size at least ${2|V(H)|\over 5}$.
\end{abstract}

\section{Introduction}

A subset $S$ of vertices in a hypergraph $H$ is {\em independent} if there are no edges of $H$ inside $S$. The cardinality of a largest independent set  of $H$ is denoted by $\alpha(H)$. A {\em linear cycle} (often also called loose cycle) in a hypergraph is a sequence of at least three edges where only the cyclically consecutive edges intersect and they intersect in exactly one vertex. Our original motivation was to prove the following conjecture that is still open.

\begin{conjecture}[Gy\'arf\'as-G.N.S\'ark\"ozy, \cite{GyarSarko}]\label{GyarSar}
One can partition the vertex set of every $3$-uniform hypergraph $H$
into $\alpha(H)$ linear cycles, edges and subsets of hyperedges.
\end{conjecture}

Note that  Conjecture \ref{GyarSar} would extend P\'osa theorem,
see \cite{LovCombEx} from graphs to $3$-uniform hypergraphs.  Conjecture \ref{GyarSar} in a weaker form (with weak cycles instead of linear cycles) has been proved in \cite{GyarSarko}. It is important that subsets of hyperedges are allowed in Conjecture \ref{GyarSar}, such an example is the complete hypergraph $K_5^3$.

Let $\rho(H)$ denote the minimum number of edges (or subsets of edges) needed to partition $V(H)$ and let $\chi(H)$ denote the chromatic number of $H$, the minimum number of colors in a vertex coloring of $H$ without monochromatic edges. The following result proves that Conjecture \ref{GyarSar} is true if there are no linear cycles in $H$.

\begin{theorem}\label{speccase} If $H$ is a $3$-uniform hypergraph without linear cycles, then $\rho(H)\le \alpha(H)$. Moreover,  
$\chi(H)\le 3$.
\end{theorem}

We find the family of hypergraphs without linear cycles intriguing and the purpose of this paper is to prove further results about it.

Let $H=(V,E)$  be a $3$-uniform hypergraph, for $v\in V$ the link of
$v$ in $H$ is the graph with vertex set $V$ and edge set $\{(x,y):
(v,x,y,) \in E \}$. The {\em strong degree} $d^+(v)$ for $v\in V$ is
the maximum number of independent edges in the link of $v$.

Our main results are as follows.

\begin{theorem}\label{min3} Suppose that $H$ is a $3$-uniform hypergraph with $d^+(v)\ge 3$ for all $v\in V$. Then $H$ contains a linear cycle.
\end{theorem}

Theorem \ref{min3} can be easily strengthened.

\begin{theorem}\label{min3S} Suppose that $H$ is a $3$-uniform hypergraph with $d^+(v)\ge 3$ for all but at most one $v\in V$. Then $H$ contains a linear cycle.
\end{theorem}

Indeed, if a graph $G$ is a counterexample with exceptional vertex $w$ to Theorem \ref{min3S} then three copies of $G$ can be joined together through cut point $w$ to get a counterexample to Theorem \ref{min3} as well.
Notice that Theorem \ref{min3S} does not hold with {\em two} exceptional vertices:  for odd $n$ consider cyclically consecutive triples of $[n]$ together with two vertices $x,y$ and with edges $xyi$ for all $i\in [n]$. This hypergraph has no linear cycles and $d^+(i)=3,d^+(x)=d^+(y)=2$.

It is worth mentioning that the condition $d^+(v)\ge 3$ cannot be weakened by requiring that the link of $v$ cannot be pierced by at most two vertices. Indeed, $K_5^3$ or hypergraphs obtained by attaching further $K_5^3$s to it are examples. It is also interesting to note that the maximal number of edges in a $3$-uniform hypergraph without linear cycles is ${n-1 \choose 2}$, the maximum number of edges without a linear triangle  \cite{CSK}, \cite{FF}.

\begin{theorem}\label{alfa2/5}
If $\cH_n$ is an $n$-vertex hypergraph without linear cycles, then $\alpha(\cH_n)\ge {2\over 5}n$.
\end{theorem}

The hypergraph consisting of vertex disjoint copies of $K_5^3$ shows that equality can hold in Theorem \ref{alfa2/5}.

\subsection{Skeletons, near-skeletons}

A {\em linear tree} is a $3$-uniform hypergraph that is obtained
from a single edge by repeatedly adding edges that intersect the previous hypergraph in exactly one vertex. A single vertex is a {\em
trivial tree}.  A {\em linear path} is a linear tree built so that the next edge always intersects the previous edge in a vertex of degree one. A {\em linear cycle} is obtained from a linear path of at least two edges, by adding an edge that intersects the first and the last edges of the path in one of their degree one vertices. For brevity, we often just use the term tree for a linear tree.

The {\em star} of a tree $T$ at $v\in V(T)$ is the subtree of $T$
containing the edges of $T$ incident to $v$. Considering the pairs
covered by the edges of $T$ as a graph $G(T)$, for any $v\in V(T)$
the pairs $(x,y)$ that are at equal distance from $v$ in $G(T)$ are
called pairs opposite to $v$. Clearly, every edge of $T$ has exactly one pair
opposite to $v$.

A {\em skeleton} $T$  in $H$ is a non-trivial subtree which cannot be extended to a larger subtree
by adding an edge $e\in E(H)$ for which $|e\cap V(T)|=1$. A  {\em near-skeleton} $T$
with an exceptional vertex $v\in V(T)$ is a {\em non-trivial} subtree
$T$ with the following property: if $|e\cap V(T)|=1$ for some $e\in
E(H)$ then $e\cap V(T)=\{v\}$. Note that skeletons are not necessarily maximum subtrees, for example
in the hypergraph with edge set $\{abc,bcd,cde\}$,  $\{bcd\}$ and  $\{abc,cde\}$ are both skeletons.
The following easy lemma is stated without proof.

\begin{lemma}\label{skelprop} Suppose $H$ is a $3$-uniform
hypergraph having no linear cycle and $T$ is a linear subtree in it. Let $v\in V(T)$ and
$f=(v,a,b)\in E(H)$ be such that  $\{a,b\}$ intersects $V(T)$ but does
not intersect the star at $v\in V(T)$. Then $\{a,b\}$ is a pair
opposite to $v$ in $T$. Replacing the edge of $T$ containing $a,b$
by $f$ is called a swap, it gives another linear tree on vertex set
$V(T)$.
\end{lemma}

The following is a useful corollary of Lemma \ref{skelprop}.

\begin{corollary} \label{span1}Suppose $T$ is a skeleton (near-skeleton) in a $3$-uniform
hypergraph $H$ that has no linear cycle. Then any sequence of swaps
with edges of $E(H[V(T)])$ results in a skeleton (near-skeleton)
$T'$ in $H$ with $V(T')=V(T)$.
\end{corollary}

\subsection{Proof of Theorem \ref{speccase}}
Consider a $3$-uniform hypergraph $H$ and choose a skeleton $T_1$ in it, then let $T_2$ be a skeleton in $H\setminus T_1$ and continue
with $T_3,\dots,T_m$ until an edgeless $T_{m+1}$ remains. Let $G_i$ be the graph obtained from $T_i$ by replacing each edge of $T_i$ by
three pairs. Observe that $\alpha(G_i)=\theta(G_i)$ where $\theta(G)$ is the minimum number of complete subgraphs whose vertices cover
$V(G)$. By the definition of skeletons, no edge of $H$ intersects $V(T_i)$ in one vertex and intersects $V(H)\setminus (\cup_{j\le i}
V(T_i))$ in two vertices.  Suppose $S_i\subset V(G_i)$ is an independent set of $G_i$. Because $H$ has no linear cycles, no edge of $H$
is inside $S_i$ and no edge of $H$ contains two vertices of $S_i$ and one vertex of $V(H)\setminus S_i$. Thus
$$\alpha(H)\ge \alpha(\cup_{i=1}^{m+1} G_i)=\sum_{i=1}^{m+1} \alpha(G_i)=\sum_{i=1}^{m+1} \theta(G_i)\ge \sum_{i=1}^{m+1} \rho(T_i)\ge \rho(H)$$  proving the first part of Theorem \ref{speccase}. The second part, $\chi(H)\le 3$, follows from $\chi(G_i)=3$ for $1\le i \le m$ and $\chi(G_{m+1})=1$, using the remarks above, that union of independent sets of $G_i$s are independent in $H$.  In fact, one can also derive $\chi(H)\le 3$ by induction, since Theorem \ref{min3} ensures that there is a vertex of $H$ with strong degree at most two. \qed

\section{Proof of Theorem \ref{min3}}

We shall prove Theorem \ref{min3} in the following slightly stronger form.

\begin{theorem}\label{ervin} Suppose that $T$ is a near-skeleton in a
$3$-uniform hypergraph $H$ and $d^+_H(v)\ge 3$ holds for every $v\in
V(T)$. Then $H$ contains a linear
 cycle.
\end{theorem}

%Note that if in a $3$-uniform hypergraph $H$ with $d^+(v)\ge 3$ for all $v\in V(H)$, we can select a near-skeleton and obtain Theorem \ref{min3}.

\noindent {\bf Proof. }
 Consider a minimum counterexample where $|V(H)|$ is as small as possible
and within that $|V(T)|$ is as small as possible. The subhypergraph of $H$ with vertex set $V(T)$ is denoted by $H(T)$.
\BalAbra{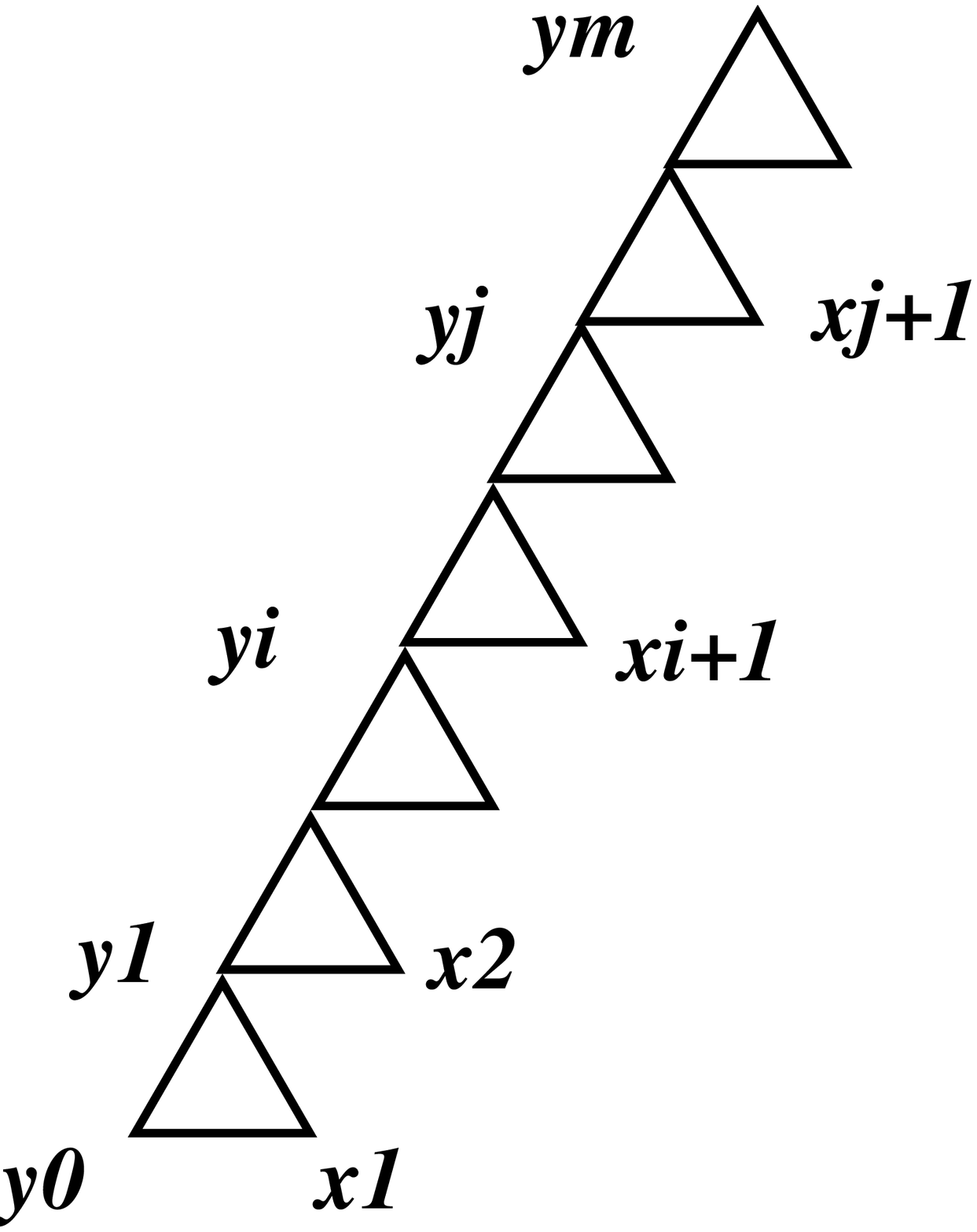}{30}{Path}{PathFig}
We may suppose that $T$ has the longest linear path $P$ among all
near-skeletons $T'$ of $H$ with $V(T')=V(T)$. Set
\begin{eqnarray*}P=\{e_1=(y_0,x_1,y_1),e_2=(y_1,x_2,y_2)\dots,e_m=(y_{m-1},x_m,y_m)\}.\end{eqnarray*}
We can see $P$ on Figure \ref{PathFig}.
%The pair $(y_{i-1},x_i)$ is called the {\em base pair} of $e_i$.
By the symmetry of $y_0,x_1$ in $P$ we may assume that $x_1$ is not the exceptional vertex of $T$.
For $1\le i<j\le m$ an {\em upward path} $B$ from $e_i$ to $e_j$ is
a linear path in $H(T)$ whose {\em first edge} intersects $e_i$ in
$\{x_i\}$, its {\em last edge} intersects $e_j$ in the pair
$\{x_j,y_j\}$ and its other vertices (inner vertices) are not on $P$.
It is possible that $B$ is a one edge path $(x_i,x_j,y_j)\in
E(H(T))$, in this case it is considered as a last edge (with no
inner vertices). A set of upward paths are {\em internally
disjoint} if their sets of inner vertices are pairwise disjoint.

\begin{figure}[ht]
\begin{center}
  \epsfig{file=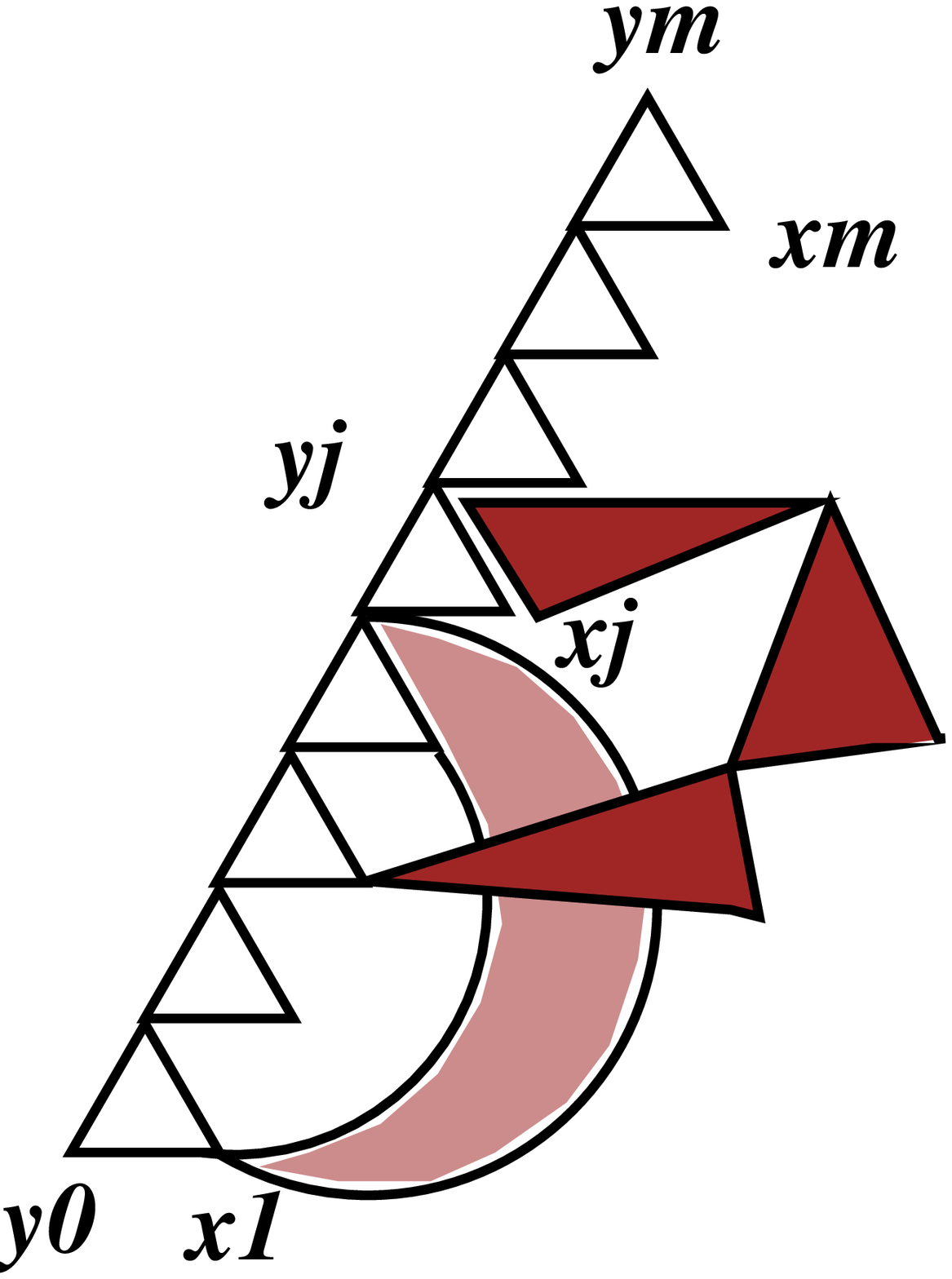,height=4cm,width=4cm}
\end{center}
\caption{Ladder}
\label{SimpleLadderB}
\end{figure}

\begin{definition}
 For $2\le j\le m$ a {\em ladder }$L_j$  is the subhypergraph of
$H(T)$ containing the path $e_1,\dots,e_j$ and a set of internally
disjoint upward paths with the following property.

\begin{itemize}
\item For every $1\le i < j$ there exists an upward path from $e_k$ to $e_{\ell}$
for some $k,\ell$ such that $1\le k \le i < \ell \le j$.
\end{itemize}
\end{definition}
Figure \ref{SimpleLadderB} shows a ladder with two upward paths.
We shall use the ladder to ensure that for any vertex $q$ not on the ladder the edge $(q,y_{i-1}y_i)$ can be continued to get a simple path from the edges of the ladder ending with a last edge of an upward path in the pair $(x_j,y_j)$.

\begin{figure}[ht]
\begin{center}
  \epsfig{file=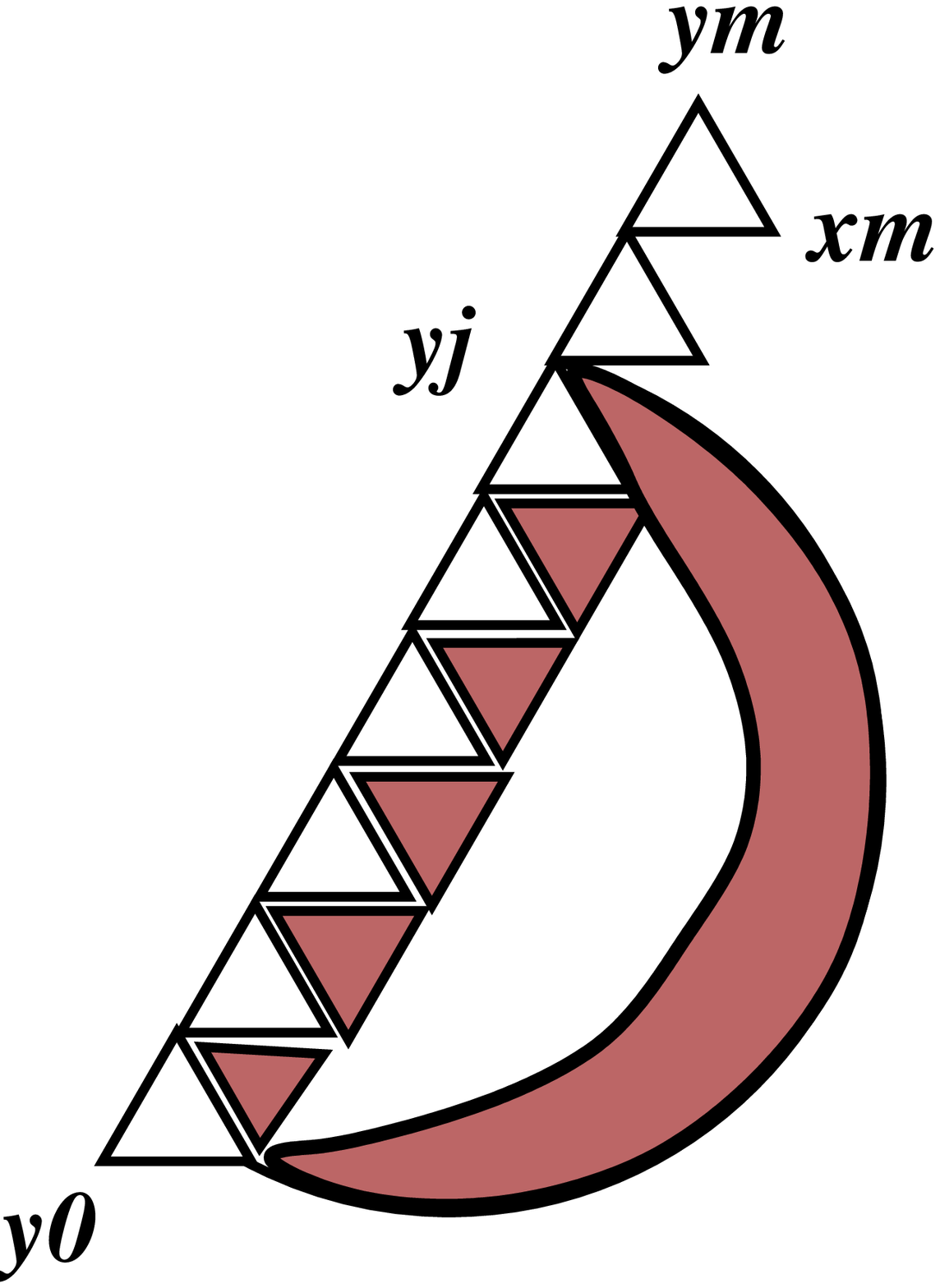,height=36mm,width=4cm}
\end{center}
\caption{A simple ladder with edges forced by Claim \ref{stepdown}}
\label{SimpleLadderC}
\end{figure}

Observe that by removing from $L_j$ the last edges of its upward
paths, we have a linear tree in $H(T)$ denoted by $L_j^*$. Ladders exist
because $d^+(x_1)\ge 3$ implies that there is an edge $f=(x_1,a,b)$ in
$H(T)$ for which $\{a,b\}\cap \{y_0,y_1\}=\emptyset$. The choice of $x_1$ and
Lemma \ref{skelprop} implies that $\{a,b\}=\{x_j,y_j\}$ for some $2\le j\le m$. Thus $P\cup f$
is a ladder $L_j$, see Figure \ref{SimpleLadderC}.

%\epsfig{file=GyFig0003.eps,width=6cm}
%\includegraphics{GyFig0003.pdf}

%Note that Claim \ref{minskel} concludes the inductive proof of Theorem \ref{ervin}
%unless if $V(T)=V(H)$ and $j=m$.  But this case can be treated as follows. Since $x_1,x_m$ are not
%exceptional, their strong degrees are at least three and that easily
%leads to contradiction ($(x_1,y_0,y_2),(x_m,x_1,y_0)\in E(H[T])$ are
%forced for example).

%\begin{claim}\label{minskel}There exist a near-skeleton $F \in H[V(T)]$ such that $V(F)$ is a proper subset of $V(T)$.
%\end{claim}

\bigskip

%\noindent \bf Proof of Claim \ref{minskel}. \rm

Let $L_j$ be a ladder such that $j$ is as large as possible. Set $P'=\cup_{i>j} e_i$ and let $M$ denote the linear
tree $P'\cup L_j^*$. We extend $M$ to a larger tree by adding a maximal linear subtree $F=F(x_j)$ of $H(T)$ with root
$x_j$,  so that its vertices (except its root) is
in $V(T)\setminus V(M)$. Notice that from the construction,
$U=V(M)\cup V(F)\subseteq V(T)$ and $M\cup F$ is a linear tree. (One can define $F$ step by step using Corollary \ref{span1}.)

Let $q\in V(F)$ and suppose that there exists $h=(q,a,b)\in E(H)$ such that $\{a,b\}\cap
V(F)=\emptyset$. The maximality of $F$ implies that $\{a,b\}\cap
M\ne \emptyset$. Applying Lemma \ref{skelprop} to the linear tree in $M\cup F$ at
vertex $q$, we get that $\{a,b\}$ either intersects the star at $q$
or it is a pair opposite to $q$. We have the following possibilities
for $\{a,b\}$.

\begin{itemize}\itemsep=0pt\parskip=0pt
\item Case 1. $\{a,b\}=\{x_k,y_k\}$ for some $k>j$
\item Case 2. $\{a,b\}=\{y_{j-1},y_j\}$
\item Case 3. Either $\{a,b\}=\{y_{k-1},x_k\}$ with some
$1\le k <j$ or $\{a,b\}$ is on an upward path of $L_j$
\end{itemize}

Case 1 would contradict to the choice of $j$ since the path with
first edge starting at $x_j$ and last edge $(q,a,b)$ would be an
upward path extending  the ladder $L_j$ to a ladder $L_k$.

Cases 2,3  for $q\ne x_j$ are also impossible since we could get a
linear cycle from the definition of the ladder $L_j$. Indeed, in
Case 2  one can start with $h=(q,a,b)$ and descend on $P$ until an
upward path leads back directly or through a jump on $P$ to
$(x_j,y_j)$, closing a cycle at $\{y_{j-1},y_j\}$. In Case 3 one can
proceed similarly but upon reaching $(x_j,y_j)$ get back to $q\in
V(F)$, closing the cycle.

We conclude that there is no $q\in V(F(x_j))\setminus \{x_j\}$ and $h=(q,a,b)\in E(H(T))$ such that $\{a,b\}\cap V(F)=\emptyset$.
Thus, if $F(x_j)\ne\{x_j\}$, $F(x_j)$ is a near-skeleton with exceptional vertex $x_j$, contradicting to the assumption.

If $F(x_j)=\{x_j\}$, the assumption $d^+(x_j)\ge 3$ allows to select $h=(x_j,a,b)\in E(H(T))$ such that $\{a,b\}\cap
\{y_{j-1},y_j\}=\emptyset$. Then $\{a,b\}$ must satisfy Case 3 and we get a linear cycle and a contradiction except when
$h=(x_j,y_0,x_1)$ and $L_j$ consists of only one upward path with one edge $f=(x_1,x_j,y_j)$ because in this case the cycle starting
with edge $(x_j,y_{k-1},x_k)$ and ending with edge $(x_1,x_j,y_j)$ degenerates. From here we assume that $L_j$ is this simple ladder
shown on Figure \ref{SimpleLadderC}.

%\begin{itemize}
%\item Fact 1. $h=(x_j,y_0,x_1)\in E(H), f=(x_1,x_j,y_j)\in E(H)$.
%\item Fact 2. for $2\le k \le j$ and $a\notin P$, $e=(a,y_{k-1},x_k)\notin E(H)$.
%\item Fact 3. for $3\le k \le j$, $e'=(a,y_{k-2},x_k)\notin E(H)$.
%\item Fact 4. $j\ge 3$.
%\end{itemize}

%Facts 2 and 3 are true for $k=j$ otherwise $e,e_{j-1},\dots,e_1,f$ or $e',e_{j-2},\dots,e_1,f$ would be a linear cycle.

%\noindent {\bf Proof of Facts 2, 3. \rm} If not true, we get a linear cycle either with $e_{k-1},e_{k-1},\dots,e_1$

In case of $j=2$ the link of $x_2$ consists of the $\{a,b\}$ pairs that are either pairs of $e_1$ or intersect $y_2$ because if $\{a,b\}=\{u,y_1\}$
with $u\notin V(P)$ then $(u,y_1,x_2),e_1,f$ would form a linear triangle. Thus, from $d^+(x_2)\ge 3$, there is an edge of $H(T)$ on $x_2$ that is
different from $h$ and does not intersect  $\{y_1,y_2\}$ and therefore would extend $L_2$ to a higher ladder. Thus we have $j\ge 3$.

\BalAbra{GyFig2}{26}{Ladder}{SimpleLadderB}
For $2\le i\le j$ define a maximal subtree $F(x_i)$ of $H(T)$ with root
$x_i$,  such that its vertices (except its root) are
in $V(T)\setminus V(M)$.

\begin{claim}\label{stepdown} For $2\le i\le j$, $F(x_i)=\{x_i\},g_i=(x_i,x_{i-1},y_{i-1})\in E(H)$ and for $3\le i \le j$,
$(x_{i-1},x_1,y_0)\notin E(H)$.
\end{claim}

\silent

\bigskip

\noindent \bf Proof of Claim \ref{stepdown}. \rm For $i=j$, $F(x_j)=\{x_j\}$. Note that for $a\notin P$, $e=(a,y_{j-1},x_j)\notin E(H)$
and $e'=(y_{j-1},y_{j-2},x_j)\notin E(H)$, otherwise $e,e_{j-1},\dots,e_1,f$ or $e',e_{j-2},\dots,e_1,f$ would be a linear cycle.

Using this and  $d^+(x_j)\ge 3$, it follows that $g_j\in E(H)$. Then $(x_{j-1},x_1,y_0)\notin E(H)$, otherwise that edge with $g_j,f$ would form a linear triangle.
This proves the claim for $i=j$. Suppose that the claim is true for some $i\ge 3$, we show it remains true for $i-1$ as well.

Suppose $F=F(x_{i-1})\ne \{x_{i-1}\}$. Then, as before, $F$ is a near-skeleton with exceptional vertex $x_{i-1}$,  a contradiction.

\BalAbra{GyFig3}{28}{A simple ladder with edges forced by Claim~\ref{stepdown}}{SimpleLadderC}

Indeed, assuming that there exists $q\in V(F)\setminus \{x_{i-1}\}$ and $h=(q,a,b)\in E(H)$ such that $\{a,b\}\cap
V(F)=\emptyset$ for some $q\in V(F)$, from Lemma \ref{skelprop} we get the following possibilities
for $\{a,b\}$.

\silent

\begin{itemize}\itemsep=0pt\parskip=0pt
\item Case A. $\{a,b\}=\{x_k,y_k\}$ for some $k>i-1$
\item Case B. $\{a,b\}=\{y_{i-2},y_{i-1}\}$
\item Case C. $\{a,b\}=\{y_{k-1},x_k\}$ with some
$1\le k \le i-1$
\end{itemize}

Case A would contradict to the choice of $j$ if $k>j$: the path with
first edge starting at $x_i$ and last edge $(q,a,b)$ would be an
upward path extending  the ladder $L_j$ to a ladder $L_k$.  If $i\le k \le j$ then  the linear cycle
starting with the path $h,g_j,\dots,g_i$ and ending with the linear path of $F$ from $x_{i-1}$ to $q$, a contradiction.

Cases B,C  are also impossible since we could get a
linear cycle. Indeed, in both cases one can start with $h=(q,a,b)$ and go up on $g_i,\dots,g_j$, return on $h$ and close the cycle
on $e_1,\dots,e_{i-1}$. Thus $F$ is a near-skeleton with exceptional vertex $x_{i-1}$ leading to contradiction.
Therefore $F(x_{i-1})=\{x_{i-1}\}$ as claimed. Moreover, $(x_{i-1},x_1,y_0)\notin E(H)$ otherwise that edge with $g_i\dots,g_j,f$ would form a linear cycle.

Now we use $d^+(x_{i-1})\ge 3$. Since $F_{i-1}=\{x_{i-1}\}$, every edge $(x_{i-1},a,b)\in E(H)$ intersects $V(P)$ and from Lemma \ref{skelprop}, we have Cases A,B,C plus those where the star at $x_{i-1}$ intersects $\{a,b\}$, i.e. $\{a,b\}\cap \{y_{i-2},y_{i-1}\}\ne \emptyset$. Notice that for $a\notin P$, $e=(a,y_{i-2},x_{i-1})\notin E(H)$ and $e'=(y_{i-2},y_{i-3},x_{i-1})\notin E(H)$
otherwise $e,e_{i-2},\dots,e_1,f,g_j,\dots,g_i$ or $e',e_{i-3},\dots,e_1,f,g_j,\dots,g_i$ would be a linear cycle. One can easily check that only three cases left for which there is no linear cycle: $\{a,b\}\cap \{y_{i-1}\}\ne \emptyset$, or $\{a,b\}=\{x_i,y_i\}$, or $\{a,b\}=\{y_{i-2},x_2\}$. From $d^+(x_{i-1})\ge 3$, all of these possibilities must occur, in particular the third, so $g_{i-1}\in E(H)$ and this completes the proof of Claim \ref{stepdown}. \qed

Observing that $g_j,\dots,g_2,f$ is a linear cycle, the proof of Theorem \ref{ervin} is completed. \qed

\section{Proof of Theorem \ref{alfa2/5}}

Let $H$ be a 3-uniform hypergraph of $n$ vertices not containing any linear cycle.
We prove that $\alpha(H) \geq 2n/5$

To facilitate the constructive proof, a {\em mixed tree}  is defined as an extension of linear 3-uniform trees where we
allow 2-element edges as well. In particular,  graph trees and 3-uniform
(linear) trees are both mixed  linear trees.  A mixed forest is the vertex-disjoint union
of mixed trees.

{\em A path-ending} of a mixed forest $T$ is a path with two edges $g,h$ where $h$ is a pendant edge (i.e. the vertices in $h\setminus g$ are of degree one in $T$) and the vertex $g\cap h$ has degree 2 in $T$. There are 4 types of path endings, depending whether $g$ or $h$ has 2 or 3 vertices. In fact, we define a degenerate
path-ending as well: in the one-edge tree the only edge is considered as a path-ending.

{\em A star-ending} of a mixed forest is a set of at least two pendant edges with a common vertex.
We state the following obvious lemma without proof.
\begin{lemma} Any mixed forest with at least one edge has either a path-ending or a star ending.
\end{lemma}

Let $T_1$ be a  maximum nontrivial skeleton, i.e. a linear tree in $H$ such that $|V(T_1)|$ is maximum. We prove the theorem constructing (step by step, details are in Subsection \ref{sconstruction})  an independent set $S$ of $H$ and a set $Z\subset V(H)$ such that

$$S\cap Z=\emptyset,\qquad {|S|\over |S|+|Z|}\ge {2\over 5},\mbox{   and   }S\cup Z=V(\cH).$$

%$$S\cap Z=\emptyset, {|S|\over |S|+|Z|}\ge {2\over 5},S\cup Z=V(H).$$

Initially set $S=Z=\emptyset$. First we cover $V(T_1)$ with $S\cup Z$ in several steps (see  Subsection \ref{sconstruction}) so that $S\subset V(T_1)$ is an independent set in $H$. Then we iterate the process, taking a maximum skeleton $T_2$ in the subhypergraph of $H$ induced by  $X=V(H)\setminus (S\cup Z)$ and continue with $T_3,T_4,\dots,T_m$, etc. until the subhypergraph of $H$ induced by $X=V(H)\setminus (S\cup Z)$ has no edges. At this point $S$ is extended by $X$ and the construction is completed.

%The
%remaining set $V(H) - V(T_1)$ of vertices is denoted by $X$.  Then every edge $e$ of $H$ containing a vertex of $V(T_1)$
%intersects  $X$ in at most one vertex and in case of $|e\cap X|=1$ the pair $e\cap V(T)$ must be part of some edge of $T_1$.

%that will be a subset of complement of the final independent set by putting the vertices into these sets. First we take a maximal skeleton $T_1$ in H and distribute its vertices in several steps.
Suppose that we have already defined $S\cup Z$ and $T_i$ so that $S$ is independent in $H$ and $S\subset \cup_{j<i} V(T_j)$. We extend $S\cup Z$ by steps in Subsection \ref{sconstruction}, in each step using a mixed forest $T$ in $T_i$, initially $T=T_i$.
We choose a vertex set $R$ (typically, but not always a subset of $V(T)$) such that a subset $R_0 \subseteq R \bigcap V(T)$ with $|R_0| \geq 2|R|/5$ vertices will be put into the independent set $S$ and $R-R_0$ is placed in $Z$. We proceed in this way until all vertices of $T_i$ are covered by $S\cup Z$.

Note that the procedure defining $S\cup Z$ ensures the properties ${|S|\over |S|+|Z|}\ge {2\over 5},S\cup Z=V(H)$ because at each step $|R_0| \geq 2|R|/5$ and in the final step $S$ is increased by $|X|$ but $Z$ is left untouched. Thus we need to ensure only that the final $S$ is independent. This will be done in Claim \ref{endclaim}.

\subsection{Construction of $S$ and $Z$}\label{sconstruction}

{\bf Case 1.} $T\subset T_i$ has a path-ending $Q=g\cup h$ with $|h|=3$. Set $h=abc$, $g=bde$ (or $g=bd$ if $|g|=2$ or $g=\emptyset$ if $T$ has one edge $h$).

{\bf Case 1.1.} There is no edge $abx$ or $bcx$ in $H$ for which $ x \notin \{Z\cup \{c,d,e\}\}$.  Put $a,b$ into $S$ and $c,d,e$ into $Z$ (ratio at least $2/5$). Replace $T$ by the mixed forest obtained from $T$ by deleting the vertices of $Q$.

%Then let $R=\{a,b,c\}\cup V(g)\}$.
%Suppose  that there is a hyperedge $f$ of $H$ containing $a$ or $b$ but not meeting $\{c,d,e\}$. By definition of Case 1.1.1, there is no such $f$ containing both $a$ and $b$.
%Suppose that $f=axy$ (or similarly $f=bxy$) with $x,y \epsilon V(H)-R$. If $x,y \epsilon  V(H)-V(T)$ then we can add $axy$ ($bxy$) to $T$, a contradiction to its maximality. If $\{x,y\} \bigcap V(T) \neq \emptyset$ then we are done by outer triple rule unless $x,y \epsilon g$ for some $g \epsilon E(T)$. We will take care of this later.

If Case 1.1 does not hold, we must have edges $abx_1,bcx_2$ in $H$ such that $ x_1,x_2 \notin \{Z\cup \{c,d,e\}\}$. However, $x_i\in V(T_i)\setminus Q$ would create a cycle in $T_i$, $x_i\in V(T_j)$ for $j<i$ would contradict the maximality of $T_j$. Thus $x_1,x_2$ are both in $X=V(H)\setminus (S\cup Z)$. If $x_1\ne x_2$ then replacing $abc$ by $abx_1$ and $bcx_2$, a skeleton larger than $T_i$ could be defined, a contradiction. Thus $x_1=x_2=x$ and we have the following.

%Then $x_i$ is not in $S\bigcap V(T)$ in because of the neighborhood rule, and not in $V(T)-S$ because it would give a cycle with the $x_i - b$ path of the skeleton (outer triple rule?). But then

{\bf Case 1.2.} There is $x\in X$ such that $abx$ and $bcx$ are edges in $H$. Put $a,c$ into $S$ and $b,x$ into $Z$ (ratio is $1/2$).  Replace $T$ by the mixed forest obtained from $T$ by deleting $\{a,b,c\}$.

%Note that if $x \in V(T)$ then it is not adjacent to the triple $abc$.
%Furthermore, by the neighborhood rule, even if $x \in S$ it is not adjacent to the triple $abc$. In both cases, take the skeleton path (of at least two edges) from $x$ to $\{a.b,c\}$ and extending it with either $abx$ or $bcx$ we get a linear cycle, a contradiction.

%So, we may assume that $x \in X$.

{\bf Case 2.} $T$ has a path-ending $Q=g\cup h$ with $|h|=2$, set $ab=h$. Put the degree one vertex of $h$ into $S$
and the other vertex of $h$ into $Z$ (ratio is $1/2$).  Replace $T$ by the mixed forest obtained from
$T$ by deleting the vertices of $h$.

{\bf Case 3.} $T$ has a star-ending.  Put one vertex of degree one from each edge  of the star
into $S$ (there are at least two) and put its other vertices into $Z$. (Clearly at least
2/5  of the vertices of the star go to $S$.)  Replace $T$ by the mixed
forest obtained from $T$ by deleting the vertices of the star-ending.

{\bf Case 4.} $T$ has only isolated vertices.   Place all vertices into $S$.

\begin{figure}[ht]
\begin{center}
  \epsfig{file=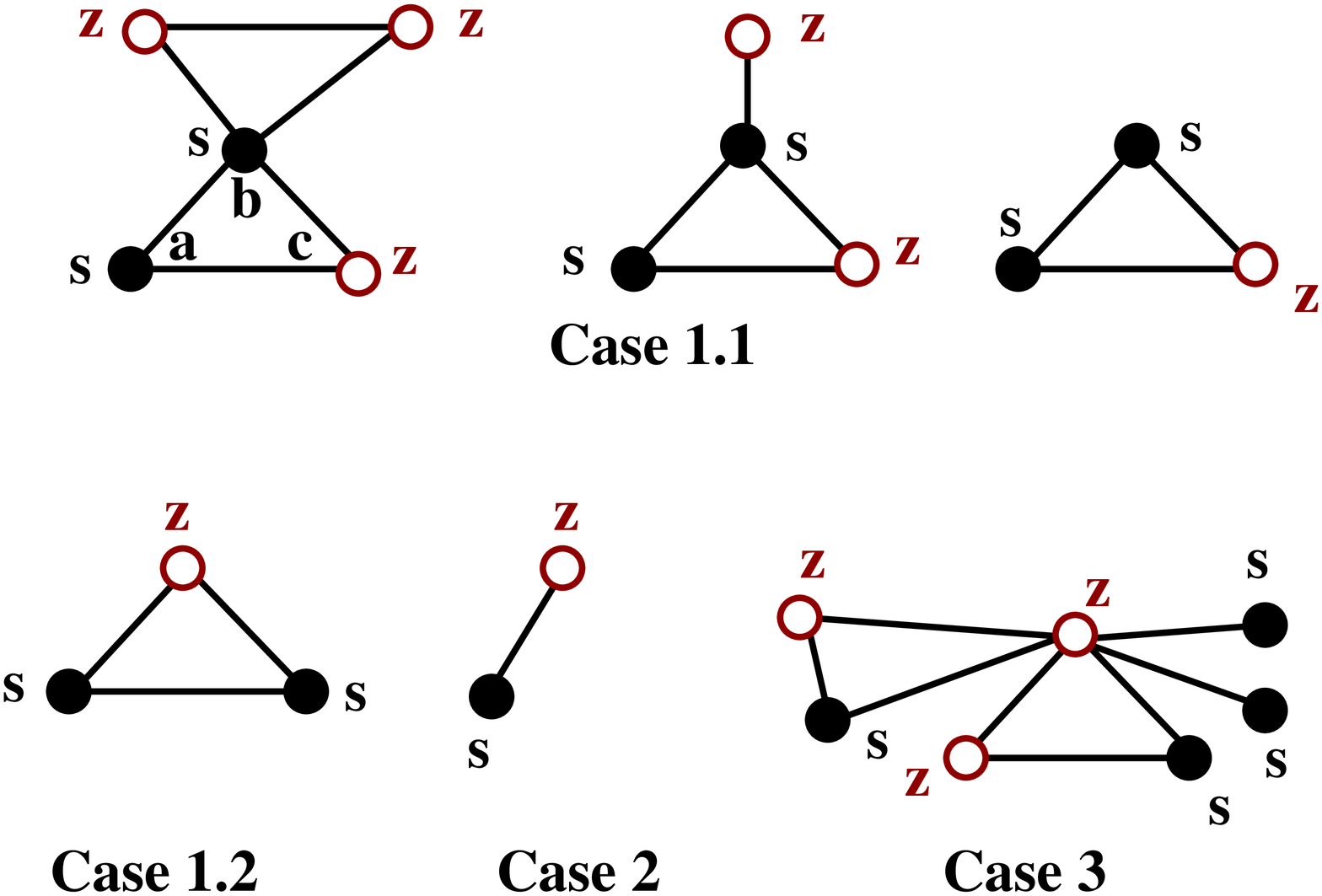,height=45mm}
\end{center}
\caption{The cases}
\end{figure}

The proof of Theorem \ref{alfa2/5} is complete with the following claim.

\begin{claim}\label{endclaim} $S$ is an independent set in $H$.
\end{claim}

\noindent {\bf Proof.}
Suppose that $e=\{s_1,s_2,s_3\}\subset S, e\in E(H)$. Observing that the construction ensures $S\subset \cup_{j=1}^m V(T_j)\cup X_m$, the maximal choices of the $T_j$s imply that the smallest $j$ for which $e\cap V(T_j)\ne \emptyset$ contains at least two vertices of $e$, say $s_1,s_2\in T_j$ and $s_3\in T_k$ for $j\le k$ or $s_3\in X_m$. If $s_3\notin V(T_j)$ then $s_3$ was placed in $S$ after $s_1,s_2$. We may assume that $s_3$ entered $S$ not earlier than $s_1,s_2$ and $s_2$ entered $S$ not earlier than $s_1$.

The $T$-neighbors of a vertex $v\in V(T)$ are the vertices in the edges of $T$ containing $v$. Notice that in Cases 2,3 the $T$-neighbors of the vertices placed in $S$ are all placed in $Z$ and in Cases 1.1 and 1.2 the $T$-neighbors of the pair placed in $S$ are placed in $Z$ (in Case 1.1 $c,d,e$, in Case 1.2 $a,x$) - we refer to this as the {\em neighbor rule.}

If $(s_1,s_2)$ or $(s_2,s_3)$ were placed in $S$ together as $(a,b)$ in Case 1.1 then the definition of Case 1 excludes $e\in E(H)$.

Suppose that $(s_1,s_2)$ or $(s_2,s_3)$ were placed in $S$ together as $(a,c)$ in Case 1.2, let $y$ denote $s_3$ or $s_1$, depending on which pair is $(a,c)$.  Then $y \in V(H)-Z-\{a,b,c,x\}$.  If $y \in X_j$ (in this case $(a,c)=(s_1,s_2),y=s_3$),  replacing the triple $abc$ by $abx$ and $acy$ in $T_j$, we get a contradiction to the maximality of $T_j$. If $y \in S\cap V(T_i)$ and the skeleton path from $y$ to $\{a,b,c\}$ reaches $b$ first, then extending it with $abx$ and $acy$ we get a linear cycle, contradiction. If the skeleton path reaches $a$ or $c$ first, then extending it by $acy$, we get a linear cycle.

Thus we may assume that no pair of the vertices of $e$ are entered into $S$ through Cases 1.2 or 1.2, i.e. they entered $S$  either in separate steps or some of them together in Case 3. Using the neighbor rule and Lemma \ref{skelprop} with $(v,a,b)=(s_1,s_2,s_3)$, $e$ would create a cycle in $H$ and this contradiction completes the proof of the claim and the proof of Theorem \ref{alfa2/5}. \qed

\section{Conjectures, open problems}

\begin{problem}Let $H$ be a $3$-uniform hypergraph with no linear cycles with no subhypergraph $K_5^3$.
Is it true that $\chi(H)\le 2$?
\end{problem}

\begin{problem}
Can one describe the structure of 3-uniform hypergraphs with no linear cycles?
\end{problem}

\begin{problem}
Is there a stability in Theorem \ref{alfa2/5}: excluding $K_5^3$ from $H$, is $\alpha(H)\ge ({2\over 5}+c)|V(H)|$ for some constant $c>0$?
\end{problem}

\begin{problem}
Which results extend to $r$-uniform hypergraphs?
\end{problem}

\noindent {\bf Acknowledgement.}  We thank conversations with Sasha Kostochka.

\end{document}